\documentclass[12pt,draft]{amsart}
\usepackage[all]{xy}
\usepackage{amsfonts}
\usepackage{amssymb}

\newtheorem{teo}{Theorem}[section]
\newtheorem{lem}[teo]{Lemma}

\theoremstyle{definition}
\newtheorem{dfn}[teo]{Definition}
\newtheorem{rk}[teo]{Remark}

\def\<{\langle}
\def\>{\rangle}
\def\ss{\subset}

\def\a{\alpha}
\def\b{\beta}

\def\d{\delta}
\def\e{\varepsilon}
\def\l{{\lambda}}
\def\r{\rho}

\def\t{\tau}

\def\f{{\varphi}}

\def\G{{\Gamma}}

\def\C{{\mathbb C}}

\def\Z{{\mathbb Z}}

\def\A{{\mathcal A}}

\def\End{\mathop{\rm End}\nolimits}
\def\Ker{\mathop{\rm Ker}\nolimits}

\def\Prim{\operatorname{Prim}}

\def\supp{\operatorname{supp}}

\def\1{\mathbf 1}

\newcommand{\ov}[1]{\overline{#1}}

\newcommand{\wh}[1]{\widehat{#1}}

\def\Fix{\operatorname{Fix}}
\def\Rep{\operatorname{Rep}}

\def\N{{\mathbb N}}

\begin{document}

\title[Noncommutative Riesz theorem and Burnside theorem]
{Noncommutative Riesz theorem
and weak Burnside type theorem on twisted conjugacy}

\author{Evgenij Troitsky}
\thanks{The author is partially supported by
RFFI Grant 05-01-00923, Grant for the support of
leading scientific schools $\rm HI\! I\! I$-619.2003.1
and Grant ``Universities of Russia'' YP.04.02.530.}
\address{Dept. of Mech. and Math., Moscow State University,
119992 GSP-2  Moscow, Russia}
\email{troitsky@mech.math.msu.su}
\urladdr{
http://mech.math.msu.su/\~{}troitsky}

\begin{abstract}
The present paper consists of two parts. In the first part, we prove
a noncommutative analogue of the Riesz(-Markov-Kakutani) theorem
on representation of functionals on an algebra of continuous
functions by regular measures on the underlying space.

In the second part, using this result, we prove a weak version of
Burnside type theorem for twisted conjugacy for arbitrary discrete groups.
\end{abstract}

\maketitle


\section{Introduction and statement of results}

\begin{dfn}
Let $G$ be a countable discrete group, and let $\phi: G\rightarrow G$ be an
endomorphism.
Two elements $x,x'\in G$ are said to be
 $\phi$-{\em conjugate} (or {\em twisted conjugate})
iff there exists a $g \in G$ such that
$$
x'=g  x   \phi(g^{-1}).
$$
We shall write $\{x\}_\phi$ for the $\phi$-{\em conjugacy} or
{\em twisted conjugacy} class
 of an element $x\in G$.
The number of $\phi$-conjugacy classes is called the {\em Reidemeister number}
of $\phi$ and is  denoted by $R(\phi)$.
If $\phi$ is the identity map, then the $\phi$-conjugacy classes are the usual
conjugacy classes in $G$.
\end{dfn}

If $G$ is a finite group, then the classical Burnside theorem (e.g., see
\cite[p.~140]{Kirillov})
says that the number of conjugacy
classes of elements of $G$ is equal to the number of equivalence
classes of irreducible representations, i.e., points of  the {\em unitary dual} $\wh G$.

Consider an automorphism $\phi$ of a finite group $G$. Then $R(\phi)$
is equal to the dimension of the space of twisted invariant functions on
$G$. Hence, by Peter-Weyl theorem (which asserts the existence of
a two-side equivariant isomorphism
$C^*(G)\cong \bigoplus_{\r\in\wh G} \End(H_\r)$),
$R(\phi)$ coincides
with the sum of dimensions
$d_\r$ of the spaces of twisted invariant elements in $\End(H_\r)$, where $\r$ runs over
$\wh G$ and the space of a representation $\r$ is denoted by $H_\r$.
By the Schur lemma,
$d_\r=1$, if $\r$ is a fixed point of $\wh\phi:\wh G\to \wh G$, where
$\wh\phi (\r):= \r\circ\phi$, and is zero otherwise. Hence,
$R(\phi)$ coincides with the number of fixed points of $\wh\phi$
(see, e.g., \cite{FelshB}). The purpose of the present paper is to generalize this
statement to the case of infinite discrete groups (after an appropriate adaptation).

\begin{rk}
If $\phi: G\to G$ is an epimorphism, then it induces the map $\wh\phi:\wh G\to\wh G$,
$\wh\phi (\r)=\r\circ\phi$
(since a representation is irreducible if and only if the  scalar operators in
the representation space
are the only ones which commute with all
representation operators). This is not the case for a general endomorphism $\phi$,
since $\r \phi$ may be reducible for an irreducible representation $\r$,
and  $\wh\phi$ can be defined only as a multi-valued map.
Nevertheless, we can define the set $\Fix (\wh\phi)$ of fixed points
of $\wh\phi$ on $\wh G$ for a general endomorphism $\phi$.
\end{rk}

\begin{dfn}
Let $\Rep(G)$ be the space of equivalence classes of
finite dimensional unitary representations of $G$.
Then the corresponding map $\wh\phi_R:\Rep(G)\to \Rep(G)$
is defined in the same way as above: $\wh\phi_R(\r)=\r\circ\phi$.

Let $\Fix(\wh\phi)$ be the set of points $\r\in\wh G\ss
\Rep(G)$ such that $\wh\phi_R (\r)=\r$.
\end{dfn}

\begin{teo}[twisted Burnside theorem for type I groups\ \cite{FelTro}]
\label{teo:mainth1} Suppose that $G$ is a finitely generated discrete group
of type {\rm I}, $\phi$ is an endomorphism of $G$, $R(\phi)$ is the
number of $\phi$-conjugacy classes, and $S(\phi)=\# \Fix
(\wh\phi)$ is the number of $\wh\phi$-invariant equivalence classes
of irreducible unitary representations. If  $R(\phi)$ or
$S(\phi)$ is finite, then $R(\phi)=S(\phi)$.
\end{teo}

Let $\mu(d)$, $d\in\N$, be the {\em M\"obius function};
i.e.,
$$
\mu(d) =
\left\{
\begin{array}{ll}
1 & {\rm if}\ d=1,  \\
(-1)^k & {\rm if}\ d\ {\rm is\ a\ product\ of}\ k\ {\rm distinct\ primes,}\\
0 & {\rm if}\ d\ \mbox{is not square-free.}
\end{array}
\right.
$$

\begin{teo}[congruences for the Reidemeister numbers \cite{FelTro}]\label{teo:mainth3}
Let $\phi:G$ $\to G$ be an endomorphism of a countable discrete group $G$
such that
all numbers $R(\phi^n)$ are finite, and let $H\ss G$ be a subgroup such that
$
  \phi(H) \subset H,
$
and for each $x\in G$ there exists an $n\in \N$  with $\phi^n(x)\in H$.
If the pair  $(H,\phi^n)$ satisfies the assumptions of
Theorem~{\rm~\ref{teo:mainth1}}
for each $n\in\N$,
then
 $$
 \sum_{d\mid n} \mu(d)\cdot R(\phi^{n/d}) \equiv 0 \mod n
 $$
for all $n$.
\end{teo}

The situation is much more complicated
for groups of type II${}_1$. For example, for
the semi-direct product defined by the action
of $\Z$ on $\Z\oplus \Z$ by a hyperbolic automorphism, we have
found  an automorphism
with finite
Reidemeister number (equal to four) such that $\wh\phi$ has at least five fixed
points on $\wh G$ \cite{FelTroVer}.
This
phenomenon is due to bad separation properties of $\wh G$ for general
discrete groups. A deeper study leads to the following general theorem, which is
one of two main results of the present paper.

 \begin{teo}[{weak Burnside type theorem for twisted classes}]
The dimension $R_*(\phi)$ of the space of
twisted invariant functions on $G$ lying in the Fourier-Stieltjes
algebra $B(G)$ is  equal to the number $S_*(\phi)$ of generalized
fixed points $I$ of the homeomorphism $\wh\phi$
$($the sum of codimensions of the subspaces generated by elements
of the form $a-\d_g * a * \d_{\phi(g^{-1})}+I$ $)$
on the Glimm spectrum of $G$, i.e., on
the complete regularization of $\wh G$, provided that at least one of
the numbers $R_*(\phi)$ and
$S_*(\phi)$ is finite. Here $\d_g$ is the delta function supported at $g$.
 \end{teo}
 This result
allows one to obtain the strong form $R(\phi)=S(\phi)$ of the twisted
Burnside theorem in a number of cases.
The proof of the generalized Burnside theorem
in~\cite{FelTro}
for groups of type I (see Theorem \ref{teo:mainth1} in the present paper
and also~\cite{FelshB,FelTroVer})
used an identification of $R(\phi)$ with the dimension of the
space of  twisted invariant ($L^\infty$-)functions on $G$, i.e., twisted
invariant functionals on $L^1(G)$. Since only part of $L^\infty$-functions
 (namely, Fourier-Stieltjes functions) define functionals on $C^*(G)$,
one a priori has
$R_*(\phi)\le R(\phi)$. Nevertheless, functions satisfying some
symmetry conditions very often lie in the Fourier-Stieltjes algebra,
so one can conjecture that $R(\phi)=R_*(\phi)$ provided that $R(\phi)<\infty$.
This is the case for all known examples.

The weak theorem will be proved as follows.
The well-known Riesz(-Markov-Kakutani) theorem identifies
the space of linear functionals on algebra $A=C(X)$ with the
space of regular measures on $X$. To prove the weak twisted Burnside
theorem, we first obtain a generalization of the Riesz theorem to the
case of a noncommutative $C^*$-algebra $A$ using the Dauns-Hofmann
theorem on representation of $C^*$-algebras by sections.
The corresponding measures
on the Glimm spectrum are functional-valued.
In the extreme situations (a commutative algebra or an algebra with one-point Glimm spectrum),
this theorem
(which is the second main result of the present paper) either is
reduced to the Riesz theorem or becomes tautological, but for
the group $C^*$-algebras of discrete groups one in many cases
obtains a new tool for counting twisted conjugacy classes.

\medskip
The interest in twisted conjugacy relations has its origins, in particular,
in the Nielsen-Reidemeister fixed point theory (e.g., see \cite{Jiang,FelshB}),
Selberg theory (e.g., see \cite{Shokra,Arthur}),
and  algebraic geometry (e.g., see \cite{Groth}).
Note that the Reidemeister number of an
endomorphism of a finitely generated Abelian group is known to be
finite if and only if $1$ is not in the
spectrum of the restriction of this endomorphism to the free part of the group
(e.g., see \cite{Jiang}).  The Reidemeister number
of any automorphism of a  nonelementary
Gromov hyperbolic group is infinite
\cite{FelPOMI}.

\medskip\noindent
{\bf Acknowledgement.}
The present research is part of our joint research programm
with A.~Fel'shtyn in Max-Planck-Institut f\"ur Mathemetik (MPI)
in Bonn.
The author thanks the MPI for kind support and
hospitality during the time in which most of this research has been completed.
I am also indebted to MPI and the organizers of the
Workshops on Noncommutative Geometry and Number Theory I and II (Bonn,
August 2003 and June 2004),
where considerable part of the results of this paper were presented.

The author is grateful to
B.~M.~Bekka,
A.~L.~Fel'shtyn,
V.~M.~Ma\-nui\-lov,
A.~S.~Mish\-che\-n\-ko,
A.~I.~Shtern,
L.~I.~Vainerman, and
A.~M.~Vershik
for helpful discussions and to the referee for valuable remarks.

\section{Algebras of operator fields}\label{sec:algfiel}

First, we recall some facts from  the theory of operator fields following
\cite{Fell} (see also \cite[\S 10]{DixmierEng}). Let $T$ be a topological space,
and let a $C^*$-algebra (or, more generally, involutive Banach algebra)
$A_t$ be assigned to  each point $t\in T$.

\begin{dfn}
A {\em continuity structure for $T$ and the family $\{A_t\}$\/} is a linear
space $F$
of operator fields on $T$ ranging in $\{A_t\}$,
(i.e., maps sending each $t\in T$ to an element of $A_t$)
and possessing the following properties:
\begin{enumerate}
\item if $x\in F$, then the real-valued function $t\mapsto \|x(t)\|$ is continuous
on $T$;
\item the set $\{x(t)\,|\,x\in F\}$ is dense in $A_t$ for each $t\in T$;
\item $F$ is closed under pointwise multiplication and involution.
\end{enumerate}
\end{dfn}

\begin{dfn}
An operator field $x$ is said to be {\em continuous\/} with respect to $F$ at
a point $t_0$
if for each $\e>0$ there exists an element $y\in F$ and a neighborhood $U$
of $t_0$ such that $\|x(t)-y(t)\|<\e$ for all $t\in U$. The field $x$ is
{\em continuous on\/} $T$ if it is continuous at all points of $T$.
\end{dfn}

\begin{dfn}
A {\em full algebra of operator fields\/} is a family $A$ of operator
fields on $T$ satisfying the following conditions
\begin{enumerate}
\item $A$ is a *-algebra, i.e., is closed under all the pointwise
algebraic operations;
\item for each $x\in A$, the function $t\mapsto \|x(t)\|$ is continuous
on $T$ and vanishes at infinity;
\item for each $t$, $\{x(t)\; | \; x\in A\}$ is dense in $A_t$;
\item $A$ is complete in the norm $\|x\|=\sup\limits_t\|x(t)\|$.
\end{enumerate}

\end{dfn}

A full algebra of operator fields is obviously a continuity
structure. If $F$ is any continuity structure, let us define
$C_0(F)$ to be the family of all operator fields $x$
continuous on $T$ with respect to $F$ such that $t\mapsto
\|x\|$ vanishes at infinity. One can prove that
$C_0(F)$ is a full algebra of operator fields.

\begin{lem}
For any full algebra $A$ of operator fields on $T$, the following three
conditions are equivalent:
\begin{enumerate}
\item $A$ is a maximal full algebra of operator fields;
\item $A=C_0(F)$ for some continuity structure $F$;
\item $A=C_0(A)$.
\end{enumerate}
\end{lem}

Such a maximal full algebra $A$ of operator fields may sometimes be called a
{\em continuous direct sum\/} of the family $\{A_t\}$.
We shall study the unital case; thus $T$ is compact, and the
property of vanishing at infinity is void. Moreover, we suppose that $T$ is Hausdorff and
hence normal. Clearly, in this case the
full algebra is
{\em separating\/}  in the sense that  if $s,t\in T$, $s\ne t$, $\a\in A_s$, and
$\b\in A_t$, then there exists an $x\in A$ such that $x(s)=\a$ and $x(t)=\b$.
We will need to distinguish an algebra $A$ itself from its
realization as the algebra $\G (\A)$ of section of the field
 $\A=\{A_t\}$ of algebras.
We denote the section corresponding to an element
$a\in A$ by $\wh a$.

\section{Functionals and measures}\label{sec:funstmeas}

\begin{dfn}\label{dfn:aofm}
Let $\Sigma$ be an algebra of subsets of $T$.
A {\em measure associated with a maximal full algebra $A=\Gamma(\A)$
of operator fields\/}
is a set function
$\mu:S\mapsto \mu(S) \in \Gamma(\A)^*=A^*$, $S\in \Sigma$,
such that
$\mu(S)(a)=0$ whenever $\supp \wh a \cap S=\emptyset$.
In what follows, we use the abbreviation AOFM for such a measure.

An AOFM $\mu$ is \emph{additive} if
$\mu(\sqcup S_i)(a)=\sum_i \mu(S_i)(a)$

An AOFM $\mu$ is \emph{bounded} if the supremum $\|\mu\|$ of
$\sum_i \|\mu(S_i)\|$ over partitions $\{S_i\}$ of $T$ is finite.

A bounded additive AOFM will be abbreviated as BA AOFM.
\end{dfn}

\begin{dfn}\label{dfn:regul}
A BA AOFM is *-{\em weakly regular\/} (a RBA AOFM)
if for any $E\in\Sigma$,
$a\in A$,
and $\e>0$ there exists a
set $F\in\Sigma$ whose closure is contained in $E$ and a set $G\in\Sigma$
whose interior contains $E$ such that $|\mu(C)a|<\e$ for every
$C\in \Sigma$ with $C\ss G\setminus F$.
\end{dfn}

In what follows, for $\Sigma$ we take the algebra of all subsets
of $T$ or the algebra generated by closed subsets of $T$.

\begin{dfn}\label{dfn:lambmn}
Suppose that an AOFM $\l$ is defined on an algebra $\Sigma$ of subsets of $T$
and $\l(\emptyset)=0$. A set $E\in\Sigma$
is called $\l$-{\em set\/}
if
$$
\l(M)=\l(M\cap E)+\l(M\cap(T\setminus E))
$$
for each $M\in\Sigma$.
\end{dfn}

\begin{lem}\label{lem:mumn}
Suppose that $\l$ is an {\rm AOFM} defined on an algebra $\Sigma$ of subsets of
$T$ and $\l(\emptyset)=0$. The family of $\l$-sets is a subalgebra
of $\Sigma$ on which $\l$ is additive.
Furthermore, if $E$ is the union of a finite set $\{E_n\}$ of
disjoint $\l$-sets and $M\in \Sigma$, then
$\l(M\cap E)=\sum_n \l(M\cap E_n).$
\end{lem}

\begin{proof}
It is clear that the void set, the whole space, and the
complement of any $\l$-set are $\l$-sets. Now let $X$ and $Y$ be
$\l$-sets, and let $M\in\Sigma$. Since $X$ is a $\l$-set, we have
\begin{equation}\label{eq:i}
\l(M\cap Y)=\l(M\cap Y \cap X)+\l(M\cap Y\cap (T\setminus X)),
\end{equation}
and since $Y$ is a $\l$-set, we have
\begin{equation}\label{eq:ii}
\l(M)=\l(M\cap Y)+\l(M\cap (T\setminus Y)),
\end{equation}
$$
\l(M\cap (T\setminus (X\cap Y)))=
\l(M\cap (T\setminus (X\cap Y))\cap Y)+
\l(M\cap (T\setminus (X\cap Y))\cap (T\setminus Y));
$$
hence
\begin{equation}\label{eq:iii}
\l(M\cap (T\setminus (X\cap Y)))=
\l(M\cap (T\setminus X)\cap Y)+
\l(M\cap (T\setminus Y)).
\end{equation}
It follows from (\ref{eq:i}) and (\ref{eq:ii}) that
$$
\l(M)=\l(M\cap Y \cap X)+\l(M\cap Y\cap (T\setminus X))
+\l(M\cap (T\setminus Y)),
$$
and (\ref{eq:iii}) implies that
$$
\l(M)=\l(M\cap Y \cap X)+\l(M\cap (T\setminus (X\cap Y))).
$$
Thus $X\cap Y$ is a $\l$-set. Since
$\cup X_n=T\setminus \cap (T\setminus X_n)$, we conclude that
the $\l$-sets form an algebra. Now if $E_1$ and $E_2$ are disjoint
$\l$-sets, then replacing $M$ by $M\cap(E_1\cup E_2)$
in Definition \ref{dfn:lambmn}, we see that
$$
\l(M\cap(E_1\cup E_2))=\l(M\cap E_1)+\l(M\cap E_2).
$$
The second assertion of the lemma follows from this by induction.
\end{proof}

It is well known that each functional $\t$ on a $C^*$-algebra $B$
can be represented as a linear combination of four positive functionals
in the following canonical way. First, represent
$\t$ as $\t=\t_1+ i\t_2$, where the self-adjoint functionals
$\t_1$ and $\t_2$ are given by
\begin{equation}\label{eq:razlnasamosop}
\t_1(a)=\frac{\t(a)+\ov{\t(a^*)}}2,\qquad \t_2(a)=\frac{\t(a)-\ov{\t(a^*)}}{2i}.
\end{equation}
By the lemma on the Jordan decomposition, each self-adjoint functional
$\a$ can be uniquely represented as the difference $\a=\a_+-\a_-$ of two positive
functionals such that
\begin{equation}\label{eq:razlozhjord}
    \|\a\|=\|\a_+\|+\|\a_-\|
\end{equation}
(see \cite[\S 3.3]{Murphy},
\cite[Theorem 3.2.5]{Ped}). Let us decompose an AOFM accordingly.
Since the decomposition is unique, it follows that the summands are AOFM.
If we start from BA AOFM, then the additivity of the summands will follow from
the uniqueness of the decomposition, and the boundedness (with a constant
twice as large)
will follow from (\ref{eq:razlnasamosop}) and property
(\ref{eq:razlozhjord}).
The same argument shows that the summands are *-weak regular, provided that
so is the original AOFM.
Thus the AOFM in the decomposition are
{\em positive\/} in the sense that
$$
\mu(E)(a^*a)\ge 0
$$
for each $E\in \Sigma$.
Such a set function is nondecreasing with the respect to set inclusion.

\begin{lem}\label{lem:weakreg1}
The sets $F$ and $G$ in the {\rm Definition~\ref{dfn:regul}}
can be chosen in such a way that
$|\mu(C)(fa)|<\e$ for each continuous
function $f:T\to [0,1]$.
\end{lem}

\begin{proof}
Consider the decompositions $\mu=\sum_{i=1}^4 x_i \mu_i$ and
$a=\sum_{j=1}^4 y_j a_j$, where $\mu_i$ and $a_j$ are positive and
$x_i$ and $y_j$ are complex numbers of norm $\le 1$. Let us choose
the sets $F$ and $G$ as in Definition~\ref{dfn:regul} for
$\e/16$ and for all pairs $\mu_i$, $a_j$ simultaneously. Then
$$
0\le
\mu_i(C)(f\cdot a_j)=\mu_i(C)((a_j)^{1/2} f(a_j)^{1/2})\le
\mu_i(C)(a_j)\le \frac \e{16},
$$
and
$$
|\mu_i(C)(f\cdot a)|\le \sum_{i,j=1}^4 |x_i y_j|\cdot
|\mu_i(C)(a_j)|\le 16\cdot\frac \e{16}=\e.
$$
\end{proof}

\begin{teo}\label{teo:razryv}
Let a unital separable $C^*$-algebra $A$ be isomorphic to
a full algebra $\G(\A)$ of operator fields over a Hausdorff space $X$.
Then the functionals on $A\cong \G(\A)$ can be identified with {\rm RBA AOFM\/}
associated with $\G(\A)$.
\end{teo}

\begin{proof}
The assumptions of the theorem imply that $T$ is a
separable Hausdorff compact and the unit ball of the dual space of $A$
is a metrizable compact set in the $*$-weak topology.

Obviously, RBA AOFM form a linear normed space with respect to
$\|.\|$.

First, we wish to prove that the natural linear map $\mu\mapsto \mu(T)$ is
an isometry of the space of RBA AOFM into $A^*$. Since $\|\mu(T)\|\le \|\mu\|$, it is
of norm $\le 1$. Now take an arbitrary small $\e>0$. Let $E_1,\dots,
E_n$ be a partition of $T$ such that
$$
\sum_{i=1}^n \|\mu(E_i)\|\ge \|\mu\|-\e.
$$
Let $a_i\in A$ be elements of norm 1 such that
$\mu(E_i)(a_i)\ge \|\mu(E_i)\|-\e/n$.

By the *-weak regularity of $\mu$ and the normality of $T$, there exist closed
sets $C_i$, disjoint open sets $G_i$, and continuous functions
$f_i:T\to [0,1]$ such that $C_i \ss E_i$,
$|\mu(E_i\setminus C_i)(a_j)|\le \e/n^2$,
$C_i\ss G_i$, $|\mu(G_i\setminus C_i)(a_j)|\le \e/n^2$
(and estimations remain valid after multiplication by positive functions
as well, as in Lemma~\ref{lem:weakreg1}),
$f_i(s)=0$ if $s\not\in G_i$, and $f_i(s)=1$ if $s\in C_i$, $i,j=1,\dots,n$.

Consider the element
$a:=\sum_i f_i a_i \in \G(\A)=A$. Then $\|a\|\le 1$ and
$$
|\:\mu(S)(a)-\|\mu\|\:|\le \sum_{i=1}^n\left| \mu(E_i)(a)-\mu(E_i) \right|+\e
$$
$$
\le \sum_{i=1}^n\left| \mu(E_i\setminus C_i)(a) + \mu(C_i)(a)
-\mu(E_i)(a_i) \right|+2\e
$$
$$
= \sum_{i=1}^n\left|  \sum_{j=1}^n\mu(E_i\setminus C_i)(f_j a_j)
+ \mu(C_i)(a_i)-\mu(E_i)(a_i) \right|+2\e
$$
$$
\le \sum_{i,j=1}^n | \mu(E_i\setminus C_i)(f_j a_j)|
+ \sum_{i=1}^n \left| \mu(E_i\setminus C_i)(a_i) \right|+2\e
\le n^2 \:\frac \e{n^2} + n \:\frac \e{n^2} + 2\e
\le 4\e.
$$
Since $\e$ is arbitrary small, it follows that $\|\mu\|=\|\mu(S)\|$.

It remains to represent the general functional $\f$ as an RBA AOFM.
This functional on $\G(\A)$ can be extended by the Hahn-Banach
theorem to a continuous functional $\psi$ on $B(\A)=\prod_{t\in T}A_t$
(the $C^*$-algebra
of not necessary continuous sections of $\A$ with the $\sup$-norm).
This functional can be decomposed as
$\psi=\sum_{i=1}^4 \a_i\psi_i$, where the $\psi_i$ are
positive functionals, $\a_i\in\C$, $|\a_i|=1$, and $\|\psi_i\|\le \|\psi\|$.
Let
$$
\l(E)(a):=\psi(\chi_E a),\qquad \l_i(a):=\psi_i(\chi_E a),\quad
i=1,\dots,4,
$$
where $a\in \G(\A)$ and $\chi_E$ is the characteristic function
of $E$. Obviously, $\l(T)(a)=\psi(a)$ and $\l$ is a BA AOFM.
Indeed, the first two properties of Definition \ref{dfn:aofm} are
obvious. The third property can be verified for each $\l_i$,
$i=1,\dots,4$:
$$
\sum_{j=1}^N |\l_i(E_j)|=\sum_{j=1}^N \l_i(E_j)(\1)=
\sum_{j=1}^N \psi_i(\chi_{E_j}\1)=\psi_i(\1)\le \|\psi_i\|;
$$
hence,
$$
\sum_{j=1}^N |\l(E_j)| \le \sum_{i=1}^4 \sum_{j=1}^N
|\l_i(E_j)|\le \sum_{i=1}^4 \|\psi_i\| \le 4\cdot \|\psi\|.
$$

Now we wish to find an RBA AOFM $\mu$ such that
$\mu(T)(a)=\l(T)(a)$. Without loss of generality, it suffices
to do this for a positive measure $\l=\l_i$.

Let $F$ represent a general closed subset, $G$ a general
open subset, and $E$ a general subset of $T$. Define $\mu_1$ and
$\mu_2$ by setting
$$
\mu_1(F)(a^*a)=\inf_{G\supset F} \l(G) (a^*a),\qquad
\mu_2(E) (a^*a)=\sup_{F\ss E} \mu_1(F) (a^*a),
$$
and then by taking the linear extension.
More precisely, owing to separability there exists a cofinal
sequence $\{G_i\}$, $G_i\supset F$. The unit ball in the dual space is weakly compact,
and there exists a weakly convergent sequence $\l(G_{i_k})$.
Its limit $\psi$ is a positive functional on $A$ enjoying the desired $\inf$-property
on positive elements. In particular, it is independent of the choice
of $\{G_i\}$ and $\{G_{i_k}\}$. One defines $\mu_2$ in a similar way.

These set functions $\mu_1$ and $\mu_2$ are nonnegative and nondecreasing.
Let $G_1$ be open, and let $F_1$ be closed. If
$G\supset F_1\setminus G_1$, then $G_1\cap G\supset F_1$
and $\l(G_1)\le \l(G_1)+\l(G)$, so that
$\mu_1(F_1)\le \l(G_1)+\l(G)$. Since $G$ is an arbitrary open
set containing $F_1\setminus G_1$, we have
$$
\mu_1(F_1)\le \l(G_1) +\mu_1(F_1\setminus G_1).
$$
If $F$ is a closed set, then from this inequality, by allowing
$G_1$ to range over all open sets containing $F\cap F_1$, we have
$$
\mu_1(F_1)\le \mu_1 (F\cap F_1) + \mu_2 (F_1\setminus F).
$$
If $E$ is an arbitrary subset of $T$ and $F_1$ ranges over all closed
subsets of $E$, then it follows from the preceding inequality that
\begin{equation}\label{eq:idanf}
\mu_2(E) \le \mu_2(E\cap F) + \mu_2(E\setminus F).
\end{equation}
We claim that
\begin{equation}\label{eq:iidanf}
\mu_2(E) \ge \mu_2(E\cap F) + \mu_2(E\setminus F)
\end{equation}
for an arbitrary subset $E$ of $T$
and arbitrary closed set $F$ in $T$.
To see this, let $F_1$ and $F_2$ be disjoint closed sets. Since $T$ is
normal, it follows that there are disjoint neighborhoods $G_1$
and $G_2$ of $F_1$ and $F_2$,
respectively. If $G$ is an arbitrary neighborhood of $F_1\cup F_2$,
then $\l(G)\ge \l(G\cap G_1) + \l(G\cap G_2)$, so that
$$
\mu_1(F_1\cap F_2)\ge \mu_1 (F_1) + \mu_2(F_2).
$$
We now let $E$ and $F$ be arbitrary sets in $T$ with $F$ closed, and
let $F_1$  range over closed subsets of $E\cap F$ while $F_2$ ranges over
the closed subsets of $E\setminus F$. The preceding inequality then
proves (\ref{eq:iidanf}). From (\ref{eq:idanf}) and (\ref{eq:iidanf}),
we have
\begin{equation}\label{eq:iiidanf}
\mu_2(E)=\mu_2(E\cap F) + \mu_2(E\cap (T\setminus F))
\end{equation}
for an arbitrary $E$ in $T$ and a closed subset $F\ss T$.
The function $\mu_2$ is defined on the algebra of all subsets
of $T$, and it follows from (\ref{eq:iiidanf}) that each closed
set $F$ is a $\mu_2$-set.
If $\mu$ is the restriction of $\mu_2$ to the algebra generated
by closed sets, then it follows from Lemma \ref{lem:mumn}
that $\mu$ is additive
on this algebra. It is clear from the definition of $\mu_1$ and $\mu_2$
that $\mu_1(F)=\mu_2(F)=\mu(F)$ if $F$ is closed, and hence
$\mu(E)=\sup_{F\ss E}\mu(F)$. This shows that $\mu$ is *-weak regular,
and since $\|\mu(T)\|<\infty$, we see that $\mu$ is RBA AOFM.

Finally, by definition, $\mu(S)(a)=\l(S)(a)=\psi(a)=\f(a)$
for $a\in \G(\A)$.
\end{proof}

\section{Twisted-invariant AOFM}\label{sec:invfunsme}

The most part of the subsequent argument is valid for various representations
of algebras by operator fields, but for now we restrict ourselves to the case
of the group $C^*$-algebra of a discrete group and concentrate
ourselves on the following important representation by sections due to Dauns and
Hofmann
\cite[Corollaries 8.13, 8.14]{DaunsHofmann}.

Let $Z$ be the center of a C*-algebra $A$, and let $\wh Z$ be the space of maximal ideals
of $Z$
equipped with the standard topology. If  $I\in {\mathcal P}:=\Prim (A)$
(the space of kernels of unitary irreducible representations), then
$Z\cap I \in \wh Z$. (This follows from the fact that
the restriction to $Z$ of an irreducible representation with kernel
$I$ gives rise to a homomorphism
$Z \to \C$ and hence $Z\cap I$ is a maximal ideal.) We obtain a map
$f: {\mathcal P}\to \wh Z$. Suppose that $T:=f({\mathcal P})$.
For each $x\in T$, consider the ideal
$I_x:=\cap I$, $f(I)=x$, (the \emph{Glimm ideal}) and the field $A/I_x$ of algebras.
We have the map $a\mapsto \{x\mapsto a+I_x\}$ of the algebra
$A$ into the algebra of sections of this field. An important
result in
\cite{DaunsHofmann} is that this map is an isomorphism.
The map
$f:{\mathcal P}\to T$ is universal with respect to continuous maps  $g:{\mathcal P}\to S$
to Hausdorff spaces, i.e., any such map can be represented as
$h\circ f$ for some continuous $h:T\to S$.
The space $T$ is compact for a unital algebra.

Now we consider a countable discrete group $G$ and an
automorphism $\phi$ of $G$. Let $A=C^*(G)$. One has the twisted action
of $G$ on $A$,
$$
g[a]=\d_g * a * \d_{\phi(g^{-1})},
$$
where $\d_g$ is the delta function supported at $g$, and a similar action on functionals,
since they are realized as some functions on $G$. (It coincides with the dual action
up to the replacement $g\to g^{-1}$.)
The same action is defined on $A/I$ since the ideals are shift invariant.

\begin{dfn}\label{dfn:generReid}
The dimension $R_*(\phi)$ of the space of twisted invariant functionals
on $C^*(G)$ is called the \emph{generalized Reidemeister number of $\phi$}.
Hence $R_*(\phi)$ is the dimension of the space of twisted invariant
elements of the Fourier-Stieltijes algebra $B(G)$.

Recall that the Fourier-Stieltijes algebra of a discrete group $G$
has the following three equivalent definitions: (1) the space
of coefficients of all unitary representations of $G$ (i.e., functions
of the form  $g\mapsto \<\r(g)\xi,\eta\>$,
$\xi,\eta\in H_\r$, where $H_\r$ is the space of a unitary representation $\r$);
(2) the space of all finite linear combinations of positively definite functions
on $G$; (3) the space of bounded linear functionals on $C^*(G)$.
The commutative multiplication in these function spaces on $G$ is introduced
pointwise.
\end{dfn}

\begin{dfn}
A (Glimm) ideal $I$ is a {\em generalized fixed point of $\wh \phi$\/}
if the linear span of
elements $b-g[b]$ is not dense in $A_I=A/I$, where
$g[.]$ is the twisted action, i.e., its closure $K_I$ does not coincide
with $A_I$.
Denote by $GFP$ the set of all generalized fixed points.
\end{dfn}

If we have only finitely many such fixed
points, then the twisted invariant RBA AOFM are
concentrated in these points. Indeed, let us describe the action of $G$
on RBA AOFM in
more detail.
The action of $G$ on measures is defined by the identification of
measures with functionals on $A$.

\begin{lem}\label{lem:ograninv}
If $\mu$ corresponds to a twisted invariant functional, then
for each Borel $E\ss T$ the functional $\mu(E)$ is twisted
invariant.
\end{lem}

\begin{proof}
This is an immediate consequence of *-weak regularity.
Indeed, suppose that
$a\in A$, $g\in G$,
$\e>0$ is an arbitrary small number, and $U$ and $F$ are defined as in
Lemma~\ref{lem:weakreg1}, for $a$ and $g[a]$ simultaneously.
Take a continuous
function $f:T\to [0,1]$ with $f|_F=1$ and $f|_{T\setminus U}=0$.
Then
\begin{eqnarray*}
|\mu(E)(a-g[a])|&=&
|\mu(E\setminus F)(a)
+ \mu(F)(a)
-\mu(E\setminus F)(g[a])
- \mu(F)(g[a])|\\
&\le& |\mu(F)(a)-\mu(F)(g[a])|+2\e =|\mu(F)(fa)-\mu(F)(g[fa])|+2\e\\
&=&|\mu(U)(fa)-\mu(U\setminus F)(fa)-\mu(U)(g[fa])
+\mu(U\setminus F)(g[fa])|+2\e\\
&\le & |\mu(U)(fa)-\mu(U)(g[fa])|+4\e =|\mu(T)(fa)-\mu(T)(g[fa])|+4\e\\
& =&4\e.
\end{eqnarray*}
\end{proof}

\begin{lem}
For each twisted-invariant functional $\f$ on $C^*(G)$, the corresponding
measure $\mu$ is concentrated on the set $GFP$ of generalized fixed points.
\end{lem}

\begin{proof}
Let $\|\mu\|=1$. Suppose opposite: there exists an element $a\in A$,
$\|a\|=1$, vanishing at the generalized fixed points and satisfying $\f(a)\ne 0$.
Let $\e:=|\f(a)|>0$. In each point $t\not\in GFP$ we can find elements
$b^i_t\in A$, $g^i_t\in G$, $i=1,\dots k(t)$,
such that
$$
\|a(t)-\sum_{i=1}^{k(t)}(g^i_t[b^i_t](t)-b^i_t(t))\|<\e/4.
$$
Then there exists a
neighborhood $U_t$ of $t$ such that
$$
\|a(s)-\sum_{i=1}^{k(t)}(g^i_t[b^i_t](s)-b^i_t(s))\|<\e/2
$$
 for $s\in U_t$.  Take a finite subcover
$\{U_{t_j}\}$, $j=1,\dots,n$,
of $\{U_t\}$
and a Borel partition $E_1,\dots,E_n$ subordinated to this subcover.
Then
$$
\f(a)=\sum_{j=1}^n \mu(E_j)(a)=
\sum_{j=1}^n \mu(E_j)\left(a-\sum_{i=1}^{k(t_j)}(g^i_{t_j}[b^i_{t_j}]
-b^i_{t_j})\right)
+\sum_{j=1}^n \sum_{i=1}^{k(t_j)}\mu(E_j)(g^i_{t_j}[b^i_{t_j}]
-b^i_{t_j}).
$$
By Lemma \ref{lem:ograninv}, each summand in the second term is zero.
The absolute value of the first term does not exceed
$\sum_j\|\mu(E_j)\| \e/2 \le \|\mu\|\cdot \e/2=\e/2$.
This contradicts the fact that $|\f(a)|=\e$.
\end{proof}

Since a functional $\f$ on $A_I$ is twisted-invariant if and only if
$\Ker\f\supset K_I$, it follows that the dimension of the space of these functionals
is equal to the dimension of the space of functionals on $A_I/K_I$ and is
finite  if and only if the space $A_I/K_I$ is
finite-dimensional. In this case, the dimension of the space of
twisted invariant functionals on $A_I$ is equal to $\dim(A_I/K_I)$.

\begin{dfn}
The number
$$
S_*(\phi):=\sum_{I\in GFP} \dim (A_I/K_I)
$$
is called the \emph{generalized number $S_*(\phi)$ of fixed points} of $\wh\phi$
on the Glimm spectrum.
\end{dfn}

Since functionals associated with measures
concentrated at distinct points
are linearly independent (the space is Hausdorff),
we see that the argument above gives the following statement.
\begin{teo}[weak generalized Burnside theorem]
$$
R_*(\phi)=S_*(\phi)
$$
provided that one of these numbers is finite.
\end{teo}


\begin{thebibliography}{10}

\bibitem{Arthur}
J.~Arthur and L.~Clozel, \emph{Simple algebras, base change, and the advanced
  theory of the trace formula}, Princeton University Press, Princeton, NJ,
  1989. \MR{90m:22041}

\bibitem{DaunsHofmann}
J.~Dauns and K.~H. Hofmann, \emph{Representations of rings by continuous
  sections}, Mem. Amer. Math. Soc., vol.~83, Amer. Math. Soc., Providence, RI,
  1968.

\bibitem{DixmierEng}
J.~Dixmier, \emph{{C*}-algebras}, North-Holland, Amsterdam, 1982.

\bibitem{Fell}
J.~M.~G. Fell, \emph{The structure of algebras of operator fields}, Acta Math.
  \textbf{106} (1961), 233--280. \MR{29 \#1547}

\bibitem{FelshB}
A.~Fel'shtyn, \emph{Dynamical zeta functions, {N}ielsen theory and
  {R}eidemeister torsion}, Mem. Amer. Math. Soc. \textbf{147} (2000), no.~699,
  xii+146. \MR{2001a:37031}

\bibitem{FelTro}
A.~Fel'shtyn and E.~Troitsky, \emph{A twisted {B}urnside theorem for countable
  groups and {R}eidemeister numbers}, Noncommutative Geometry and Number Theory
  (C.~Consani and M.~Marcolli, eds.), Vieweg, Braunschweig, 2006, pp.~141--154
  (Preprint MPIM2004--65, math.RT/0606155).

\bibitem{FelTroVer}
A.~Fel'shtyn, E.~Troitsky, and A.~Vershik, \emph{Twisted {B}urnside theorem for
  type {II}${}_1$ groups: an example}, Preprint 85, Max-Planck-Institut f{\"u}r
  Mathematik, 2004, (to appear in \emph{Math. Res. Lett.}, math.RT/0606161).

\bibitem{FelPOMI}
A.~L. Fel'shtyn, \emph{The {R}eidemeister number of any automorphism of a
  {G}romov hyperbolic group is infinite}, Zap. Nauchn. Sem. S.-Peterburg.
  Otdel. Mat. Inst. Steklov. (POMI) \textbf{279} (2001), no.~6 (Geom. i
  Topol.), 229--240, 250. \MR{2002e:20081}

\bibitem{Groth}
A.~Grothendieck, \emph{Formules de {N}ielsen-{W}ecken et de {L}efschetz en
  g\'eom\'etrie alg\'ebrique}, S\'eminaire de G\'eom\'etrie Alg\'ebrique du
  {B}ois-{M}arie 1965-66. {SGA} 5, Lecture Notes in Math., vol. 569,
  Springer-Verlag, Berlin, 1977, pp.~407--441.

\bibitem{Jiang}
B.~Jiang, \emph{Lectures on {N}ielsen fixed point theory}, Contemp. Math.,
  vol.~14, Amer. Math. Soc., Providence, RI, 1983.

\bibitem{Kirillov}
A.~A. Kirillov, \emph{Elements of the theory of representations},
  Springer-Verlag, Berlin Heidelberg New York, 1976.

\bibitem{Murphy}
G.~J. Murphy, \emph{${C}^*$-algebras and operator theory}, Academic Press, San
  Diego, 1990.

\bibitem{Ped}
G.~K. Pedersen, \emph{{C*}-algebras and their automorphism groups}, Academic
  Press, London -- New York -- San Francisco, 1979.

\bibitem{Shokra}
Salahoddin Shokranian, \emph{The {S}elberg-{A}rthur trace formula}, Lecture
  Notes in Mathematics, vol. 1503, Springer-Verlag, Berlin, 1992, Based on
  lectures by James Arthur. \MR{MR1176101 (93j:11029)}

\end{thebibliography}

\def\cprime{$'$} \def\dbar{\leavevmode\hbox to 0pt{\hskip.2ex \accent"16\hss}d}
  \def\polhk#1{\setbox0=\hbox{#1}{\ooalign{\hidewidth
  \lower1.5ex\hbox{`}\hidewidth\crcr\unhbox0}}}
\providecommand{\bysame}{\leavevmode\hbox to3em{\hrulefill}\thinspace}
\providecommand{\MR}{\relax\ifhmode\unskip\space\fi MR }
\providecommand{\MRhref}[2]{%
  \href{http://www.ams.org/mathscinet-getitem?mr=#1}{#2}
}
\providecommand{\href}[2]{#2}

\end{document}